\documentclass{amsart}

\pdfoutput=1

\usepackage{amsmath}
\usepackage{amsfonts}
\usepackage{amsthm}
\usepackage{amssymb}

\newcommand{\titletext}{Toroidal embeddings and polyhedral divisors}

\usepackage[pdfauthor={Robert Vollmert},
			pdftitle={\titletext}]{hyperref}

\usepackage{graphicx}

\theoremstyle{plain}
\newtheorem{theoremcite}{Theorem}

\newtheorem*{theoremmain}{Theorem 1}
\newtheorem{lemma}{Lemma}
\theoremstyle{remark}
\newtheorem{remark}{Remark}
\newtheorem*{notation}{Notation}
\newtheorem*{example}{Example}

\begin{document}

\newcommand{\ZZ}{\mathbb{Z}}
\newcommand{\CC}{\mathbb{C}}
\newcommand{\QQ}{\mathbb{Q}}
\renewcommand{\AA}{\mathbb{A}}
\newcommand{\poldiv}{\mathcal{D}}
\newcommand{\cO}{\mathcal{O}}
\newcommand{\cA}{\mathcal{A}}
\newcommand{\til}{\widetilde}
\newcommand{\Spec}{\operatorname{Spec}}
\newcommand{\PPDiv}{\operatorname{PPDiv}}
\newcommand{\Star}{\operatorname{star}}
\newcommand{\Hom}{\operatorname{Hom}}
\newcommand{\dprime}{{\prime\prime}}
\newcommand{\eval}{\operatorname{eval}}
\newcommand{\TV}[1]{\operatorname{TV}(#1)}
\newcommand{\ratto}{\dashrightarrow}
\newcommand{\isoto}{\overset{\sim}{\to}}
\newcommand{\Pol}{\operatorname{Pol}}
\newcommand{\CaDiv}{\operatorname{CaDiv}}
\newcommand{\pos}{\operatorname{pos}}
\newcommand{\SL}{\operatorname{SL}}
\newcommand{\Mat}{\operatorname{Mat}}
\newcommand{\conv}{\operatorname{conv}}
\newcommand{\princdiv}{\operatorname{div}}
\newcommand{\ptxset}[2]{(#1,#2)}

\newcommand{\mins}[2]{{\min \langle #1, #2 \rangle}}

\title{\titletext}

\author[R. Vollmert]{Robert Vollmert}
\address{Fachbereich Mathematik und Informatik, 
Freie Universit\"at Berlin,
Arnimalle 3, 
14195 Berlin, 
Germany}
\email{vollmert@math.fu-berlin.de}


\begin{abstract}
  Given an effective action of an $(n-1)$-dimensional torus on an
  $n$-dimensional normal affine variety, Mumford constructs a
  toroidal embedding, while Altmann and Hausen give a description in
  terms of a polyhedral divisor on a curve. We compare the fan of the
  toroidal embedding with this polyhedral divisor.
\end{abstract}

\maketitle


\section*{Introduction}

Suppose $X$ is an $n$-dimensional normal affine variety over the complex
numbers with an effective action by the $(n-1)$-dimensional torus $T$.
With $T \cong (\CC^*)^{n-1}$, we associate the lattice $M \cong \ZZ^{n-1}$
of characters and the dual lattice $N = \Hom(M, \ZZ)$ of
one-parameter subgroups. The action defines the weight cone $\omega$
in $M$ generated by the degrees of semi-invariant functions on $X$ and
the dual cone $\sigma$ in $N$. Effectivity of the action translates to
the fact that $\omega$ is full-dimensional and $\sigma$ is pointed.
\begin{notation}
  A cone $\delta$ ``in'' a lattice $N$ is really a subset of the vector
  space $N_\QQ = N \otimes \QQ$. The toric variety associated with this cone
  will be denoted by $\TV\delta$.
\end{notation}

Our goal is to compare two sets of combinatorial data associated with
$X$. Mumford~\cite[Chapter 4, \S 1]{kkms} takes a rational quotient map
$p$ from $X$ to a complete nonsingular curve $C$. He defines $X^\dprime$
to be the
normalization of the graph of $p$ and shows that for certain open
subsets $U$ of $C$, we obtain a toroidal embedding $(U \times T,
X^\dprime)$. This determines a combinatorial datum, namely the 
toroidal fan $\Delta(X,U)$. It is a collection of cones in
\emph{different} lattices $\ZZ \times N$, one for each point $P \in C
\setminus U$, glued along their common face in $\ptxset{0}{N}$.

Altmann and Hausen~\cite{ppdiv} construct a divisor $\poldiv$ with polyhedral
coefficients on a nonsingular curve $Y$; this divisor determines a
$T$-variety $\til{X}$, affine over $Y$, which contracts to $X$. Here,
$\poldiv$ is of the form $\Sigma_{P \in Y} \Delta_P \otimes P$, where
the $\Delta_P$ are polyhedra in $N_\QQ$ with tail cone
$\sigma$, only finitely many nontrivial.

To compare these data, we note that the curve $Y$ is an open subset of
$C$, namely the image of the map $\pi \colon X^\dprime \to C$. In
fact, the varieties $\til{X}$ and $X^\dprime$ agree, which allows us
to describe $\Delta(X,U)$ in terms of $\poldiv$. Defining the
homogenization of a polyhedron $\Delta \subset N_\QQ$ with tail
$\sigma$ to be the cone in $\ZZ \times N$ generated by
$\ptxset{1}{\Delta}$ and $\ptxset{0}{\sigma}$, we obtain the following
result.

\newcommand{\theoremtext}{%
  The toroidal fan $\Delta(X,U)$ is equal to the fan obtained by
  gluing the homogenizations of the coefficient polyhedra $\Delta_P$
  of points $P \in Y \setminus U$ along their common face
  $\ptxset{0}{\sigma}$.
}

\begin{theoremmain}
  \theoremtext
\end{theoremmain}

In Section~\ref{sec:toroidal}, we recall relevant facts about toroidal
embeddings and summarize the construction of the embedding $(U \times
T, X^\dprime)$. Section~\ref{sec:polyhedral} contains some details
about polyhedral divisors on curves. Finally, we present the proof of
Theorem~1 in Section~\ref{sec:comparison}.

\section{Toroidal interpretation} \label{sec:toroidal}

\subsection*{Toroidal embeddings}

A toroidal embedding~\cite[Chapter 2]{kkms} is a pair $(U,X)$ of a
normal variety $X$ and an open subset $U \subset X$ such that for each
point $x \in X$, there exists a toric variety $(H,Z)$ with embedded
torus $H \subset Z$ which is locally formally isomorphic at some
point $z \in Z$ to $(U,X)$ at $x$. We will further assume that the
components $E_1,\dotsc,E_r$ of $X \setminus U$ are normal, i.e., that all
toroidal embeddings are ``without self-intersection''.

The components of the sets
  $\cap_{i \in I} E_i \setminus \cup_{i \not\in I} E_i$ 
for all subsets $I \subset \{1,\dotsc,r\}$ give a stratification of
$X$. The \emph{star} of a stratum $Y$ is defined to be the union of
strata $Z$ with $Y \subset \overline{Z}$.  Given a stratum $Y$, we
have the lattice $M_Y$ of Cartier divisors on the star of $Y$ with
support in the complement of $U$. The submonoid of effective divisors
is dual to a polyhedral cone $\sigma_Y$ in the dual lattice $N_Y$.

If $Z \subset \Star(Y)$ is a stratum, its cone $\sigma_Z$ is a face of
$\sigma_Y$. The \emph{toroidal fan} of the embedding $(U,X)$ is the union of
the cones $\sigma_Y$ glued along common faces.

\begin{remark}
  A toroidal fan differs from a conventional fan only in that it lacks a
  global embedding into a lattice.
\end{remark}

Below, we will use the fact that an \'etale morphism $(U,\Star(Y)) \to
(H,\TV{\delta})$ induces an isomorphism $\sigma_Y \isoto \delta$
of lattice cones.

\subsection*{Toroidal embeddings for torus actions}
We return to the $T$-variety $X$ and summarize Mumford's
description~\cite[Chapter 4, \S 1]{kkms}.  There is a canonically
defined rational quotient map $p \colon X \dashrightarrow C$ to a
complete nonsingular curve $C$.  Sufficiently small invariant open
sets $W \subset X$ split as $W \cong U \times T$ for some open set $U
\subset C$, where the first projection $U \times T \to U$ corresponds
to $p$. We will identify $U \times T$ with $W$.

We define $X^\prime$ to be the closure of the graph of the rational
map $p$ in $X \times C$, and $X^\dprime$ to be its normalization. The
action of $T$ on $X$ lifts to $X^\dprime$. We may consider $U \times
T$ as an open subset of $X^\dprime$; the projection to $U$ now
extends to a regular map $\pi \colon X^\dprime \to C$.

After possibly replacing $U$ by an open subset, we are in the
following situation: Let $P \in C \setminus U$ be a point in the
complement of $U$. The sets $U$, $U^\prime = U \cup \{P\}$ and
$\pi^{-1}(U^\prime)$ are affine with coordinate rings $R$, $R^\prime$
and $S$, respectively. We may regard $S$ as a subring of $R \otimes
\CC[M]$ which is generated by homogeneous elements with respect to the
$M$-grading.  Denoting by $s$ a local parameter at $P \in C$, the ring
$S$ is generated over $R^\prime$ by a finite number of monomials
$s^k\chi^u$.

The corresponding semigroup in $\ZZ \times M$ and its dual cone
$\delta_P$ in $\ZZ \times N$ define a toric variety
$Z = \TV{\delta_P}$. The monomial generators of $S$ define an \'etale map
$\pi^{-1}(U^\prime) \to Z$ which shows that the embedding $(U \times
T, \pi^{-1}(U^\prime))$ is toroidal with cone isomorphic to
$\delta_P$. By considering all points $P \in C \setminus U$, we see
that $(U \times T, X^\dprime)$ is a toroidal embedding.

\begin{theoremcite}[{\cite[Chapter 4, \S 1]{kkms}}] \label{thm:toroidal}
  The embedding $(U \times T, X^\dprime)$ is toroidal. Its fan 
  $\Delta(X,U)$ consists of the cones $\delta_P$ glued along the
  common face $\delta_P \cap \ptxset{0}{N_\QQ}$.
\end{theoremcite}

\begin{remark} \label{rem:complement}
  This common face is $\sigma \subset N_\QQ$ and corresponds to
  $\pi^{-1}(U)$, an open subset of each $\pi^{-1}(U^\prime)$.  For
  points $P$ that lie outside the image of $\pi$, we have
  $\pi^{-1}(U^\prime) = \pi^{-1}(U)$, hence the cone $\delta_P$ is
  equal to $\ptxset{0}{\sigma}$.
\end{remark}

\begin{remark} \label{rem:canonical}
  Given $U \subset C$, the constructed toroidal fan $\Delta(X,U)$ is
  independent of the choice of equivariant isomorphism $U \times T \cong
  W$. It does however depend on the choice of $U$.

  If we don't require that there be an \'etale model for the whole of
  $\pi^{-1}(U^\prime)$, we can enlarge $U$ to form a canonical
  embedding $(V \times T, \til{X})$. Here, $V$ is obtained by adding
  to any $U$ as above all points $P$ with a toric model that splits as
  $Z = \AA^1 \times F$, where $F = \TV{\sigma}$ is the generic fiber of
  $\pi$. That is, the points $P$ with cone $\delta_P$ isomorphic to
  $\sigma \times \QQ_{\ge 0}$.
\end{remark}

\begin{example}
  The affine threefold $X = \SL(2,\CC) = \CC[a,b,c,d]/(ad-bc-1)$ admits a
  two-dimensional torus action by defining
  \begin{equation*}
    (t_1,t_2) \cdot \begin{pmatrix} a & b \\ c & d \end{pmatrix} =
        \begin{pmatrix} t_1a & t_2b \\ t_2^{-1}c & t_1^{-1}d\end{pmatrix}.
  \end{equation*}
  It admits a quotient morphism $\pi \colon X \to \AA^1 = \Spec \CC[s]$
  with $s \mapsto ad$. Let $W$ be the open subset of matrices with no
  vanishing entries. With $U = \AA^1 \setminus \{0,1\}$, we get an
  isomorphism $W \cong U \times T$ by mapping $t_1 \mapsto a$ and
  $t_2 \mapsto b$.
  
  We consider $P = 0$, so $U^\prime = \AA^1 \setminus \{1\}$. The
  coordinate ring of $\pi^{-1}(U^\prime)$ is generated over
  $\CC[s]_{s(s-1)}$ by $t_1$, $st_1^{-1}$ and $t_2^{\pm1}$. Thus
  $\delta_0$ is generated by $(1,0,0)$ and $(1,1,0)$. Similarly,
  $\delta_1$ is generated by $(1,0,0)$ and $(1,0,1)$, as shown in
  Figure~\ref{fig:sl2}. The fan $\Delta(X,U)$ is obtained by gluing
  these two cones at the vertex.
    
  \begin{figure}
  \centering
  \begin{minipage}{0.3\textwidth}\centering
    \includegraphics{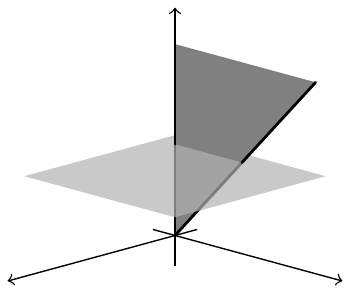}
  \end{minipage}
  \hfill
  \begin{minipage}{0.3\textwidth}\centering
    \includegraphics{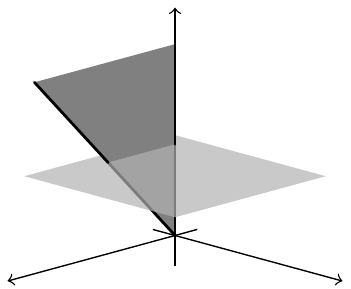}
  \end{minipage}
  \hfill
  \begin{minipage}{0.3\textwidth}\centering
  	\includegraphics{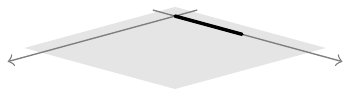} \\[3ex]
	\includegraphics{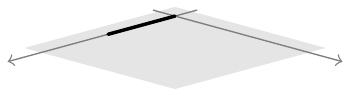}
  \end{minipage} \\

  \makebox[0.3\textwidth]{$\delta_0$} \hfill
  \makebox[0.3\textwidth]{$\delta_1$} \hfill
  \makebox[0.3\textwidth]{$\Delta_1$, $\Delta_2$}  
  \caption{Cones and polyhedra for $\SL(2,\CC)$}
  \label{fig:sl2}
  \end{figure}
\end{example}

\section{Polyhedral divisors on curves} \label{sec:polyhedral}

We turn to the construction and relevant properties of proper
polyhedral divisors on curves, restating results of Altmann and
Hausen~\cite{ppdiv} in the setting of codimension one actions.

Given a cone $\sigma$ in $N$, the set of polyhedra with tail cone
$\sigma$
\begin{equation*}
  \Pol_\sigma^+ = \{ \Delta \subset N_\QQ \mid \Delta = \Pi + \sigma 
                    \text{ for some compact polytope } \Pi \}
\end{equation*}
forms a semigroup under Minkowski addition. It is embedded in the
group of differences $\Pol_\sigma$; the neutral element is $\sigma$. A
divisor $\poldiv \in \Pol_\sigma \otimes \CaDiv(Y)$ on a smooth curve
$Y$ is called a \emph{polyhedral divisor}. Under
certain positivity assumptions ($\sum \Delta_P \subsetneq \sigma$ is
almost the right condition, see~\cite[Example 2.12]{ppdiv}), $\poldiv$
is called \emph{proper}. We may express it as
\begin{equation*}
  \poldiv = \sum \Delta_P \otimes P,
\end{equation*}
where the sum ranges over all prime divisors of $Y$, and all but
finitely many of the polyhedra $\Delta_P$ are equal to $\sigma$.

A proper polyhedral divisor defines an affine $T$-variety. Each weight
$u$ in the weight monoid $\omega \cap M$ gives a $\QQ$-divisor $\poldiv(u)$
on $Y$ by
\begin{equation*}
  \poldiv(u) = \sum \mins{u}{\Delta_P} \cdot P.
\end{equation*}
This allows us to define an $M$-graded sheaf $\cA$ of
$\cO_Y$-algebras by setting $\cA_u = \cO_Y(\poldiv(u))$.
We denote by $\til{X}$ the relative
spectrum $\Spec_Y(\cA)$ and by $X = \mathcal{X}(\poldiv)$ its
affine contraction $\Spec\Gamma(Y, \cA)$.

We summarize the relevant results on proper polyhedral divisors.
\begin{theoremcite}[{\cite[Theorem 3.4]{ppdiv}}] \label{thm:ppdiv:exist} 
  Given a $T$-variety $X$ as above, there is a curve $Y$ and a proper
  polyhedral divisor $\poldiv$ on $Y$ such that the associated
  $T$-variety $\mathcal{X}(\poldiv)$ is equivariantly isomorphic to
  $X$.
\end{theoremcite}

\begin{theoremcite}[{\cite[Theorem 3.1]{ppdiv}}] \label{thm:ppdiv:props}
  Let $X$ and $\til{X}$ be given by a proper polyhedral divisor on the
  curve $Y$.
  \begin{enumerate}
  \item The contraction map $\til{X} \to X$ is proper and birational.
  \item The map $\pi \colon \til{X} \to Y$ is a good quotient for the
    $T$-action on $\til{X}$; in particular, it is affine.
  \item \label{thm:ppdiv:affine} There is an affine open subset $U
    \subset Y$ such that the contraction map restricts to an
    isomorphism on $\pi^{-1}(U)$.
  \end{enumerate}
\end{theoremcite}

\begin{example}
  A polyhedral divisor for the torus action on $X = \SL(2,\CC)$ is computed
  easily by considering the closed embedding in the toric variety
  $\Mat(2 \times 2, \CC) \cong \AA^4$. The toric 
  computation~\cite[Section 11]{ppdiv} shows that $\AA^4$
  with the induced $(\CC^*)^2$-action may be described by the
  divisor $\poldiv^\prime = \Delta_1 \otimes D_1 
                          + \Delta_2 \otimes D_2$ on $\AA^2$, where
  $D_i = \princdiv(x_i)$ are the coordinate axes and 
  $\Delta_i = \conv\{0,e_i\}$.
  The image of $X$ in $\AA^2$ is the line through $(1,0)$ and $(0,1)$.
  Hence, $\poldiv^\prime$ restricts to the divisor 
  $\poldiv = \Delta_1 \otimes [0] + \Delta_2 \otimes [1]$
  on $\AA^1$.
\end{example}

\section{Comparison} \label{sec:comparison}

Now to compare the toroidal and polyhedral data associated with
a $T$-variety $X$. By Theorem~\ref{thm:ppdiv:exist}, 
we may assume $X$ is given by a polyhedral divisor $\poldiv$ on a
curve $Y$, contained in the complete curve $C$. As above, we have
the $T$-variety $\til{X}$ with the quotient map $\pi$ to $Y$ and
the contraction to $X$.

Denote the open subset of points $P$ with trivial coefficient
$\Delta_P = \sigma$ by $V$. Then for any open subset $U \subset V$, we
have
\begin{equation*}
  \pi^{-1}(U) = \Spec_U \cO_U \otimes \CC[\omega \cap M] 
              = U \times \TV{\sigma}.
\end{equation*}        
In particular, $U \times T$ is an open subset of
$\til{X}$. By part~\ref{thm:ppdiv:affine}
of Theorem~\ref{thm:ppdiv:props}, we may regard $U \times T$ as a
subset of $X$ after possibly shrinking $U$. The projection to $U$
gives the required rational quotient map $X \ratto C$.

We get varieties $X^\prime$ and $X^\dprime$ as before and note the
following fact.
\begin{lemma}\label{lemma:XtilXdprime}
  $\til{X}$ is canonically isomorphic to $X^\dprime$.
\end{lemma}
\begin{proof}
It follows from the construction of $X^\dprime$ that the maps
$\til{X} \to X$ and $\til{X} \to C$ factor through a map $\varphi\colon
\til{X} \to X^\dprime$. Since both maps to $X$ are proper, so is
$\varphi$. Since both maps to $C$ are affine, so is $\varphi$. Since
$\varphi$ is also birational, it is an isomorphism.
\end{proof}

Now for suitable $U$, we saw above that $(U \times T, \til{X})$ is a
toroidal embedding with fan $\Delta(X,U)$. We recall the statement
of our claim.

\begin{theoremmain}
  \theoremtext
\end{theoremmain}

To see this, consider $P \in Y \setminus U$ and $U^\prime = U \cup
\{P\}$ with local parameter $s$ at $P$. Since $\poldiv_{|U}$ is
trivial, we have $\poldiv_{|U^\prime} = \Delta_P \otimes P$. The
graded parts of $\cA = \bigoplus_{u \in \omega \cap M} \cA_u$ are
thus
\begin{equation*}
  \cA_u = \cO_{U^\prime}\big(\poldiv_{|U^\prime}(u)\big) 
        = \cO_{U^\prime}\big(\mins{u}{\Delta_P} \cdot P\big)
        = \cO_{U^\prime}\big(\lfloor \mins{u}{\Delta_P}\rfloor \cdot P\big).
\end{equation*}
Hence, we can express the graded parts of the coordinate ring $S$
of $\pi^{-1}(U^\prime)$ as
\begin{equation*}
  S_u = \Gamma\big(U^\prime, \cO_{U^\prime}(\poldiv(u))\big)
      = R^\prime \cdot s^{-\lfloor\mins{u}{\Delta_P}\rfloor}.
\end{equation*}
It follows that
the monomial semigroup of the toric model consists of the pairs
$(k,u) \in \ZZ \times M$ with $k \ge -\mins{u}{\Delta_P}$. By
Lemma~\ref{lemma:homog} below, we see that $\delta_P$ is the homogenization
of $\Delta_P$. As Remark~\ref{rem:complement} implies that points in the
complement of $Y$ don't contribute to $\Delta(X,U)$, the proof is
complete. 

\begin{lemma}\label{lemma:homog}
  Let $\Delta$ be a polyhedron in $N$ with tail cone $\sigma$. Let
  $\delta$ in $\ZZ \times N$ be its homogenization, i.e.,
  $\delta = \pos\{(0,\sigma), (1,\Delta)\}$. Then the dual cone
  $\delta^\vee$ consists of those pairs $(r,u) \in \QQ \times M_\QQ$
  with $u \in \sigma^\vee$ and $r \ge -\mins{u}{\Delta}$.
\end{lemma}

\begin{proof}
  By definition, we have $(r,u) \in \delta^\vee$ if and only if $(r,u)$
  is non-negative on both $(0,\sigma)$ and $(1,\Delta)$. The first
  condition is equivalent to $u \in \sigma^\vee$. The second condition
  means that $r \ge -\langle u,v \rangle$ for any $v \in \Delta$, that is,
  $r \ge -\mins{u}{\Delta}$.
\end{proof}

\begin{example}
  For the example of $\SL(2,\CC)$, clearly the homogenizations of the
  segments $\conv\{0,e_i\}$ give the cones $\delta_0$, $\delta_1$ generated by
  $(1,0)$ and $(1,e_i)$. This is illustrated in Figure~\ref{fig:sl2}.
\end{example}

\begin{remark}
  Both descriptions generalize to the non-affine case. Mumford treats
  this directly, while the polyhedral approach involves the fans of
  polyhedral divisors developed by Altmann, Hausen and
  S\"u\ss{}~\cite{fansydiv}. It should be straightforward to carry this
  result over.
\end{remark}

\bibliography{../bibliography/toroidal,../bibliography/polyhedral}
\bibliographystyle{hamsplain}

\end{document}